\documentclass[12pt]{amsart}
\usepackage{latexsym,amsmath,amssymb,amscd}
\usepackage{eucal}
\theoremstyle{plain}
\swapnumbers
\newtheorem{theorem}{Theorem}[section]

\newtheorem{lemma}[theorem]{Lemma}

\theoremstyle{definition}
\newtheorem{definition}[theorem]{Definition}

\theoremstyle{remark}

\newcommand{\larray}{\left(\begin{array}{cc}\right.}
\newcommand{\rarray}{\right.\end{array}\right)}

\DeclareMathOperator{\RH}{H}

\DeclareMathOperator{\HP}{HP}
\DeclareMathOperator{\End}{End}
\DeclareMathOperator{\tr}{Trace}

\DeclareMathOperator{\spin}{spin}
\newcommand{\ra}{\rightarrow}

\newcommand{\xr}{\xrightarrow}

\newcommand{\bdsplit}{\begin{displaymath}
\begin{split}}
\newcommand{\edsplit}{\end{split}
\end{displaymath}}

\newcommand{\beqn}{\begin{equation}}
\newcommand{\eeqn}{\end{equation}}
\newcommand{\bsplit}{\begin{split}}
\newcommand{\esplit}{\end{split}}

\newcommand{\bba}{\mathbb{A}}

\newcommand{\bbc}{\mathbb{C}}

\newcommand{\bbf}{\mathbb{F}}

\newcommand{\bbq}{\mathbb{Q}}

\newcommand{\bbt}{\mathbb{T}}

\newcommand{\bbv}{\mathbb{V}}

\newcommand{\bbz}{\mathbb{Z}}

\newcommand{\norm}[1]{\| #1 \| }

\newcommand{\calh}{\mathcal{H}}

\newcommand{\call}{\mathcal{L}}

\newcommand{\calo}{\mathcal{O}}

\newcommand{\cals}{\mathcal{S}}

\begin{document}
\title[Geometric Counterpart]{A geometric counterpart of the Baum-Connes
map for $GL(n)$.}
\author[J. Brodzki]{Jacek Brodzki}
\address{Faculty of Mathematical Studies \\
University of Southampton \\
Southampton SO17 1BJ\\ U.K.}
\thanks{Jacek Brodzki was supported in part by a Fellowship from the
Leverhulme Trust.}
\email{J.Brodzki@maths.soton.ac.uk}
%\author[N. Higson]{Nigel Higson}
%\address{Department of Mathematics \\
%Pennsylvania State University \\
%University Park PA 16802, U.S.A.}
%\email{higson@math.psu.edu}
\author[R. Plymen]{Roger Plymen}
\address{Department of Mathematics \\
University of Manchester \\
Manchester M13 9PL\\ U.K.}
\email{roger@maths.man.ac.uk}
\begin{abstract}
We describe a geometric counterpart of the Baum-Connes map for
the $p$-adic group $GL(n)$.
\end{abstract}
\date{\today}
\maketitle

\section*{Introduction}
Let $F$ be a nonarchimedean local field and let $G = GL(n) = GL(n,F)$.
The goal of this paper is to construct the following commutative
diagram:
$$
\begin{CD}
K_*^{\text{top}}(G)  @>{\mu}>>    K_*(C^*_r(G)) \\
@V{\text{ch}}VV                   @VV{\text{ch}}V  \\
\HP_*(\mathcal{H}(G)) @>{\imath_*}>> \HP_*(\mathcal{S}(G)) \\
@VVV     @VVV \\
\RH^*_c(\Pi(G); \mathbb{C}) @>>> \RH^*_c(\Pi^t(G);
\mathbb{C})
\end{CD}
$$
The topological $K$-theory $K^{\text{top}}_*(G)$ is defined to be
the equivariant $K$-homology of the universal example $\underline{E}G$.
As a model for $\underline{E}G$ we will take the affine building $\beta{G}$
of $G$.
Then $\beta G$  is the product of a simplicial complex by an affine line.
In addition $\beta G$ is a contractible space on which $G$ acts properly
and the quotient $\beta G/G$ is compact. Explicitly,
$$
K^{\text{top}}_*(G) = KK_*^G (C_0(\beta G), \bbc)
$$

The map $\mu$ from the topological $K$-theory of $G$ to the $K$-theory
of the reduced $C^*$-algebra of $G$ has been shown by V. Lafforgue to
be an isomorphism for all reductive $p$-adic groups \cite{L}.
The Hecke algebra $\calh (G)$ is the algebra of all complex-valued,
compactly supported, uniformly locally constant functions on $G$ with
the convolution product.
The Schwartz algebra $\cals(G)$ is the algebra of all rapid decay functions
on $G$, again with the convolution product.
There is a natural inclusion $\calh (G) \ra \cals (G)$. The periodic
cyclic homology $\HP_* (\calh (G))$ is to be understood in a purely
algebraic sense. The periodic cyclic homology
of $\HP_*(\cals (G))$ is to be understood in terms of topological
algebras and the inductive tensor products \cite{BP1}. It is proved in
\cite{BHP} that the map $\imath_*$ is an isomorphism, and a more detailed
proof is given below. The Chern character on the right hand side of the
diagram was studied in \cite{BP3} and shown to be an isomorphism
after tensoring over $\bbz$ with $\bbc$.

An outline of the construction of the Chern character on the left hand side of the diagram is
given in \cite{BHP}.   In this construction the Chern character factors through the chamber homology
of $\beta G$:

\[
ch : K^{top}_*(G) \longrightarrow \text{H}_*(G; \beta G) \longrightarrow
\text{HP}_*(\mathcal{H}(G)).
\]

We could also define the Chern character indirectly as the composite of
three maps so as to make the upper part of the diagram commutative. If
denote the Chern character on the left hand side of the diagram by
$\text{ch}_L$ and the Chern character on the right hand side of
the diagram by $\text{ch}_R$ then the definition of $\text{ch}_L$
is
$$
\text{ch}_L = \imath_*^{-1} \circ \text{ch}_R \circ \mu
$$
This creates the following commutative diagram
$$
\begin{CD}
K_*^{\text{top}}(G)\otimes_{\bbz}\bbc  @>{\mu}>>    K_*(C^*_r(G))
\otimes_{\bbz}\bbc \\
@V{\text{ch}_L}VV                   @VV{\text{ch}_R}V  \\
\HP_*(\mathcal{H}(G)) @>{\imath_*}>> \HP_*(\mathcal{S}(G)) \\
\end{CD}
$$
in which each map is an isomorphism of complex vector spaces.

From the point of view of noncommutative geometry it is natural to
seek the spaces which underlie the noncommutative algebras
$\calh(G)$ and $\cals(G)$. We prove that, at the level
of periodic cyclic homology, this spaces reveal themselves
in terms of representation theory. The smooth dual
$\Pi(G)$ has, with the the aid of Langlands parameters, a natural structure of complex manifold. The tempered dual
$\Pi^t(G)$ is, with the aid of Harish-Chandra parameters, a disjoint union of compact orbifolds. These are the two
spaces which underlie $\calh(G)$ and $\cals(G)$. Not only that
but there is a deformation retraction of the smooth dual
onto the tempered dual which, in the context of our main  commutative
diagram,  induces the Baum-Connes map.

We give a detailed description of the $q$-projection introduced in \cite{BP2}.
In the last section of this paper we track the fate of
supercuspidal representations of $G$ through the commutative diagram.
In particular, the index map $\mu$ manifests itself as an example
of Ahn reciprocity.

We would like to thank Colin Bushnell and J-F.~Dat for helpful comments,
Paul Baum and Nigel Higson for many valuable
conversations.

\section{The complex structure on the smooth dual of $GL(n)$}

The field $F$ is a nonarchimedean local field, so that
$F$ is a finite
extension of
$\bbq_p$, for some prime $p$ or $F$ is a finite
extension of the function field $\bbf _p ((x))$. The residue field $k_F$ of $F$
is the quotient $\mathbf{o}_F/\mathbf{m}_F$
of the ring of integers $\mathbf{o}_F$ by its unique maximal ideal $\mathbf{m}_F$.
Let  $q$ be the cardinality of $k_F$.

The essence of local class field theory, see \cite[p.300]{N}, is a pair
of maps $$ (d : G \longrightarrow \widehat{\mathbb{Z}}, v :
F^{\times} \longrightarrow \mathbb{Z}) $$ where $G$ is a profinite
group, $\widehat{\mathbb{Z}}$ is the profinite completion of
$\bbz$,  and $v$ is the valuation.

Let $\overline{F}$ be a separable algebraic closure of $F$. Then
the absolute Galois group
$G(\overline{F}|F)$ is the
projective limit of the finite Galois groups $G(E|F)$ taken over
the finite extensions $E$ of $F$ in $\overline{F}$. Let $\widetilde{F}$
be the maximal unramified extension of $F$.
The map $d$ is in this case the projection map
$$
d : G(\overline{F}|F) \longrightarrow G(\tilde{F}|F) \cong
\hat{\mathbb{Z}}
$$
The group $G(\tilde{F}|F)$ is procyclic.   It has a single topological
generator:  the Frobenius automorphism $\phi_F$ of $\tilde{F}|F$.
The Weil group $W_F$ is by definition the pre-image of $<\phi_F>$ in
$G(\overline{F}\mid{F})$. We thus  have the surjective map
$$
d : W_F \longrightarrow \mathbb{Z}
$$
The pre-image of $0$ is the inertia group $I_F$. In other words we have
the following short exact sequence
$$
\begin{CD}
1 @>>{}> I_F @>>{}> W_F @>>{}> \bbz @>>{}> 0
\end{CD}
$$ The group $I_F$ is given the profinite topology induced by
$G(\overline{F}|F)$.   The topology on the Weil group  $W_F$ is
dictated by the above short exact sequence.   The Weil group $W_F$
is a locally compact group with maximal compact subgroup $I_F$.
The map $$ W_F \longrightarrow G(\tilde{F}|F) $$ is a continuous
homomorphism with dense image.

A detailed account of the Weil group for local fields
may be found in \cite{T}.  For a topological group $G$ we denote
by $G^{\text{ab}}$ the quotient $G^{\text{ab}}= G/G^c$ of $G$ by
the closure $G^c$ of the commutator subgroup of $G$. Thus $G^{\text{ab}}$ is
the maximal abelian Hausdorff quotient of $G$. The local
reciprocity laws \cite[p.320]{N} $$ r_{E|F} :
G(E|F)^{ab} \cong F^{\times}/N_{E|F}E^{\times} $$ now create an
isomorphism \cite[p.69]{N2}:
$$ r_F : W^{ab}_F \cong F^{\times} $$

We have $W_F = \sqcup \Phi^{n}I_F, n \in \mathbb{Z}$.   The Weil group is a
locally compact, totally disconnected
group, whose maximal compact subgroup is $I_F$. This subgroup is also open.
There are three models for the Weil-Deligne group.

One model  is the crossed product
$W_F \ltimes \bbc$, where the Weil group acts on $\bbc$ by
$w\cdot x = \norm{w} x$, for all $w\in W_F$ and $x\in \bbc$.

The action
of $W_F$ on $\bbc$ extends to an action of $W_F$ on $SL(2,\bbc)$.
The semidirect product $W_F \ltimes SL(2,\bbc)$ is then isomorphic
to the direct product $W_F \times SL(2,\bbc)$, see
\cite[p.278]{Knapp}.   Then a complex representation of $W_F
\times SL(2,\bbc)$ is determined by its restriction to $W_F \times
SU(2)$, where $SU(2)$ is the standard compact Lie group.

From now on, we shall use this model for the Weil-Deligne group:
$$
\call_F=W_F \times SU(2).
$$

%Let $G$ be a connected reductive algebraic group over $F$. The
%$L$-group of $G$ is $$ {}^LG = {}^L G^0
%\rtimes  W_F $$ where $W_F$ is the Weil group of $F$ and
%${}^L G^0$ is the dual group of $G$. If $G'$ is
%another group over $F$ then a homomorphism $$ \rho : {}^L G \ra
%{}^L G' $$ is called an  $L$-map if $\rho$ induces the identity
%map on $W_F$, i.e. if the projection of $\rho(g , w)$ to $W_F$
%equals $w$ for any $g\in {}^L G^0$ and $w\in W_F$.

In the next definition, the complex general linear group $GL(n,
\mathbb{C})$ is equipped with the {\it discrete} topology.

\begin{definition}
 An $L$-parameter is a continuous homomorphism
$$
\phi: \call_F \ra GL(n, \mathbb{C})
$$
such that
%composition with the projection to $W_F$ induces the
%identity map on $W_F$, and
$\phi(w)$ is semisimple for all $w\in W_F$.
Two $L$-parameters are equivalent if they are conjugate under $GL(n, \mathbb{C})$.
The set of equivalence classes of $L$-parameters is denoted
$\Phi(G)$.
\end{definition}
\begin{definition}
A representation of $G$ on a complex vector space $V$ is {\it smooth} if the stabilizer of each vector
in $V$ is an open subgroup of $G$.   The set of equivalence classes of irreducible smooth
representations of $G$ is the {\it smooth dual} $\Pi(G)$ of $G$.
\end{definition}

\begin{theorem} Local Langlands Correspondence for $GL(n)$.
There is a natural bijection between $\Phi(GL(n))$ and
$\Pi(GL(n))$.
\end{theorem}

The naturality of the bijection involves compatibility of
the $L$-factors and $\epsilon$-factors attached to the two types
of objects.

%If $\rho: {}^LG \ra {}^LG'$ is an $L$ map, composition with $\rho$
%gives a transfer of $L$-parameters from $\Phi_F(G)$ to
%$\Phi_F(G')$. When the conjecture is true, this results in a
%transfer of $L$-packets on $G$ to $L$-packets on $G'$.

The local Langlands conjecture for $GL(n)$ was proved
by Stuhler \cite{Stuhler} when $F$ has positive characteristic and
by
Harris-Taylor \cite{HT} and Henniart \cite{Hen} when $F$ has
characteristic zero.
%In this case, each $L$-packet is a singleton.

We recall that a \emph{matrix coefficient} of a representation $\rho$ of a group
$G$ on a vector space $V$ is a function on $G$ of the form
$f(g) = \langle \rho(g)v, w\rangle$, where $v \in V$, $w\in V^*$, and $V^*$
denotes  the
dual space  of $V$. The inner product is given by the duality between
$V$ and $V^*$. A representation $\rho$ of $G$ is called \emph{supercuspidal}
if and only if the support of every matrix coefficient is compact modulo
the centre of $G$.

Let $\tau_j = \text{spin}(j)$ denote the $(2j + 1)$-dimensional complex
irreducible representation of the compact Lie group $SU(2)$, $j =
0, 1/2, 1, 3/2, 2, \ldots$.

For $GL(n)$ the local Langlands correspondence works in the following way.

\begin{itemize}
{\item Let $\rho$ be an irreducible representation of the Weil
group $W_F$.   Then $\pi_F(\rho \otimes 1)$ is an irreducible
supercuspidal representation of $GL(n)$, and every irreducible supercuspidal
representation of $GL(n)$ arises in this way.
 If $\det(\rho)$ is a
unitary character, then $\pi_F
(\rho \otimes 1)$ has unitary
central character, and so is pre-unitary.}

\item   We have $\pi_F(\rho \otimes \spin(j)) = Q(\Delta)$, the
Langlands quotient associated to the segment
$\{|\;|^{-(j-1)/2}\pi_F(\rho), \ldots,
|\;|^{(j-1)/2}\pi_F(\rho)\}$.   If $\det(\rho)$ is unitary, then
$Q(\Delta)$ is in the discrete series.
In particular, if $\rho = 1$ then  $\pi_F(1 \otimes \spin(j))$ is
the Steinberg representation  $St(2j + 1)$ of $GL(2j
+ 1)$.

{\item   If $\phi$ is an $L$-parameter for $GL(n)$ then $\phi =
\phi_1 \oplus \ldots \oplus \phi_m$ where $\phi_j = \rho_j \otimes
\spin(j)$.   Then $\pi_F(\rho)$ is the Langlands quotient
$Q(\Delta_1, \ldots, \Delta_m)$.   If $\det(\rho_j)$ is a unitary
character for each $j$, then $\pi_F(\phi)$ is a tempered
representation of $GL(n)$.}
\end{itemize}

This correspondence creates, as in \cite[p.~381]{Ku}, a natural
bijection
$$ \pi_F : \Phi (GL(n)) \to \Pi(GL(n)). $$

A quasi-character $\psi : W_F\ra \bbc^\times$ is {\it unramified}
if $\psi $ is trivial on the inertia group $I_F$. Recall the short
exact sequence $$ 0 \xr{} I_F \xr{} W_F \xr{d}\bbz \xr{}0 $$ Then
$\psi(w) = z^{d(w)}$ for some $z\in \bbc^\times$. Note that $\psi$
is not a \emph{Galois} representation unless $z$ has finite order in the complex torus
$\mathbb{C}^{\times}$, see \cite{T}.
Let $\Psi (W_F)$
denote the group of all unramified quasi-characters of $W_F$. Then
$$
\begin{array}{rcl}
\Psi(W_F) & \simeq & \bbc^\times \\
\psi & \mapsto & z
\end{array}
$$
Each $L$-parameter $\phi: \call_F \ra GL(n,\bbc)$ is of the form
$\phi_1 \oplus \dots \oplus \phi_m$ with each $\phi_j$
irreducible.   Each irreducible $L$-parameter is of the form $\rho
\otimes \spin(j)$ with $\rho$ an irreducible representation of the
Weil group $W_F$.

\begin{definition}\label{orbit}
The orbit $\calo(\phi) \subset \Phi_F(G)$ is defined as follows
$$ \calo(\phi) = \{ \bigoplus_{r=1}^m\psi_r\phi_r \mid \psi_r
\in \Psi(W_F), 1\leq r\leq m \}
$$
where each $\psi_r$ is an
unramified quasi-character of $W_F$.
\end{definition}

An equivalent definition of the orbit $\calo(\phi)$ :   Let
${}^LM$ be the minimal Levi subgroup of the $L$-group ${}^LG$
which contains the image $\phi(\call_F)$, see \cite[Prop. 8.6,
p.41]{Borel}. Now twist the image $\phi(\call_F)$ by all
unramified quasi-characters of ${}^LM$.   This creates the orbit
$\calo(\phi)$.

\begin{definition}
Let $\det \phi_r$ be a unitary character, $1\leq r\leq m$ and
let $\phi = \phi_1 \oplus \ldots \oplus \phi_m$. The compact orbit
$\calo^t(\phi)\subset \Phi^t(G)$ is defined as follows:
$$ \calo^t(\phi) = \{ \bigoplus_{r=1}^m\psi_r\phi_r \mid \psi_r
\in \Psi(W_F), 1\leq r\leq m \}
$$
where each $\psi_r$ is an
unramified unitary character of $W_F$.
\end{definition}

We note that $I_F \times SU(2) \subset W_F \times SU(2)$ and in
fact $I_F \times SU(2)$ is the maximal compact subgroup of
$\mathcal{L}_F$.   Now let $\phi$ be an $L$-parameter.   Moving (if
necessary) to another point in the orbit $\mathcal{O}(\phi)$ we
can write $\phi$ in the canonical form
\[
\phi = \phi_1 \oplus \ldots \oplus \phi_1 \oplus \ldots \oplus
\phi_k \oplus \ldots \oplus \phi_k
\]
where $\phi_1$ is repeated $l_1$ times, $\ldots$, $\phi_k$ is repreated $l_k$ times, and the representations
\[
\phi_j|_{I_F \times SU(2)}
\]
are irreducible and pairwise inequivalent, $1 \leq j \leq k$.   We will now write $k
= k(\phi)$.   This natural number is an invariant of the orbit
$\mathcal{O}(\phi)$.  We have
\[
\mathcal{O}(\phi) = Sym^{l_1} \mathbb{C}^{\times} \times \ldots
\times Sym^{l_k} \mathbb{C}^{\times}
\]
the product of symmetric products of $\mathbb{C}^{\times}$.

\begin{theorem}\label{torus}
The set $\Phi(GL(n))$ has the structure of complex algebraic
variety. Each irreducible component $\calo(\phi)$ is isomorphic to
the product of a complex affine space and a complex torus $$ \calo
(\phi) = \bba^l \times (\bbc^\times )^k $$
where $k = k(\phi)$.
\end{theorem}

\begin{proof}
Let $Y = \bbv(x_1y_1 - 1, \ldots, x_ny_n - 1) \subset \bbc^{2n}$.
Then $Y$ is a Zariski-closed set in $\bbc^{2n}$, and so is an
affine complex algebraic variety.   Let $X = (\bbc^{\times})^n$.
Set $\alpha : Y \rightarrow X, \alpha(x_1, y_1, \ldots, x_n,y_n) =
(x_1, \ldots, x_n)$ and $\beta : X \rightarrow Y, \beta(x_1,
\ldots, x_n) = (x_1, x_1^{-1}, \ldots, x_n, x_n^{-1})$.   So $X$
can be embedded in affine space $\bbc^{2n}$ as a Zariski-closed
subset.   Therefore $X$ is an affine algebraic variety, as in
\cite[p.50]{Sm}.

Let $A = \bbc[X]$ be the coordinate ring of $X$.  This is the
restriction to $X$ of polynomials on $\bbc^{2n}$, and so $A =
\bbc[X] = \bbc[x_1, x_1^{-1}, \ldots, x_n, x_n^{-1}]$, the ring of
Laurent polynomials in $n$ variables $x_1, \ldots, x_n$.
Let $S_n$ be the symmetric group, and let $Z$ denote the quotient variety $X/S_n$.
The variety $Z$ is an affine complex
algebraic variety.

The cooordinate ring of $Z$ is

$$\bbc[Z] \simeq \bbc[x_1, \dots , x_n, x^{-1}_1, \dots
,x^{-1}_n]^{S_n}.$$

Let  $\sigma_i$, $i= 1, \dots , n$ be the elementary symmetric
polynomials in $n$ variables.  Then from the last isomorphism we
have
\begin{displaymath}
\begin{split}
\bbc[Z] &  \simeq \bbc[x_1, \dots, x_n ] ^{S_n} \otimes \bbc
[\sigma ^{-1}_n] \\ & \simeq \bbc[\sigma _1, \dots , \sigma_n]
\otimes \bbc [\sigma^{-1}_n] \\ & \simeq \bbc [\sigma_1, \dots ,
\sigma_{n-1}]\otimes \bbc[\sigma_n , \sigma^{-1}_n] \\ & \simeq
\bbc[\bba^{n-1}]\otimes \bbc[\bba - \{0\} ]\\ & \simeq \bbc
[\bba^{n-1} \times (\bba - \{0\})]
\end{split}
\end{displaymath}

where $\mathbb{A}^n$ denotes complex affine $n$-space.  The coordinate ring of
the quotient variety $\bbc^{\times n}/S_n$ is isomorphic to the
coordinate ring of $\bba^{n - 1} \times (\bba - \{0\})$.
Now the categories of affine algebraic varieties and of finitely
generated reduced $\bbc$-algebras are equivalent, see
\cite[p.26]{Sm}.  Therefore
the variety $\bbc^{\times n}/S_n$ is isomorphic to the variety
$\bba^{n - 1} \times (\bba - \{0\})$.

Consider $\bba - \{0\} = \bbv(f)$ where $f(x) = x_1 x_2 - 1$. Then
$\partial f/\partial x_1 = x_2 \neq 0$  and $\partial f/\partial
x_2 = x_1 \neq 0$ on the variety $\bbv(f)$. So $\bba - \{0\}$ is
smooth. Then $\bba^{n-1} \times (\bba - \{0\}) $ is smooth.
Therefore the quotient variety $\bbc^{\times n}/S_n$ is a smooth
complex affine algebraic variety of dimension $n$.   Now each
orbit $\calo(\phi)$ is a product of symmetric products of
$\mathbb{C}^{\times}$.  Therefore each orbit $\calo(\phi)$ is a smooth complex affine
algebraic variety.   We have
\[
\mathcal{O}(\phi) = Sym^{l_1} \mathbb{C}^{\times} \times \ldots
\times Sym^{l_k} \mathbb{C}^{\times} = \mathbb{A}^l
\times (\mathbb{C}^{\times})^k
\]

where $l = l_1 + \ldots + l_k - k$ and $k = k(\phi)$.
\end{proof}

We now transport the complex structure from $\Phi(GL(n))$ to
$\Pi(GL(n))$ via the local Langlands correspondence.   This leads
to the next result.

\begin{theorem}
The smooth dual $\Pi(GL(n))$ has a natural complex structure.
Each irreducible component is a smooth complex affine algebraic variety.
\end{theorem}

The smooth dual $\Pi(GL(n))$ has countably many irreducible
components of each dimension $d$ with $1 \leq d \leq n$.  The
irreducible supercuspidal representations of $GL(n)$ arrange themselves into
the $1$-dimensional tori.

Each irreducible component is a smooth
affine scheme, i.e. of the form
$Spec(R)$ where $R$ is a commutative unital ring .
In fact each $R$ is a reduced finitely generated $\mathbb{C}$-algebra.
From the point of view of
noncommutative geometry, the smooth dual
$\Pi(GL(n))$ is a {\it noncommutative affine scheme} underlying
the noncommutative non-unital Hecke algebra $\calh(GL(n))$.

It follows from Theorems 1.6 and 1.7 that the smooth dual $\Pi (GL(n))$ is a complex
manifold.  Then $\bbc \times \Pi (GL(n))$ is a complex
manifold.   So the local $L$-factor $L(s, \pi_v)$ and the local
$\epsilon$-factor $\epsilon(s, \pi_v)$ are functions of {\it
several complex variables}:

$$ L : \bbc \times \Pi (GL(n)) \longrightarrow \bbc $$

$$ \epsilon : \bbc \times \Pi (GL(n)) \longrightarrow \bbc. $$

\section{Periodic cyclic homology of the Hecke algebra}

The Bernstein variety $\Omega(G)$ of $G$ is the set of $G$-conjugacy classes
of pairs $(M, \sigma)$, where $M$ is a Levi (i.e. block-diagonal) subgroup
of $G$, and $\sigma$ is an irreducible supercuspidal representation of $M$.
Each irreducible smooth representation of $G$ is a subquotient of an
induced representation $i_{GM}\sigma$. The pair $(M, \sigma)$ is unique up
to conjugacy.
This creates a finite-to-one map, the infinitesimal character,
from $\Pi(G)$ onto $\Omega G$.

Let $\Omega(G)$ be the Bernstein variety of $G$.   Each point in
$\Omega(G)$ is a conjugacy class of cuspidal pairs $(M,\sigma)$.
A quasicharacter $\psi : M \longrightarrow \mathbb{C}^{\times}$
is {\it unramified} if $\psi$ is trivial on $M^{\circ}$.   The group
of unramified quasicharacters of $M$ is denoted $\Psi(M)$.   We have
$\Psi(M) \cong (\mathbb{C}^{\times})^\ell$ where $\ell$ is the parabolic
rank of $M$.   The group $\Psi(M)$ now creates orbits:
the orbit of $(M,\sigma)$ is $\{(M, \psi \otimes \sigma): \psi \in
\Psi(M)\}$.   Denote this orbit by $D$, and set
$\Omega = D/W(M,D)$, where $W(M)$ is the Weyl group of $M$ and
$W(M,D)$ is the subgroup of $W(M)$ which leaves $D$ globally invariant.
The orbit $D$ has the structure of a complex torus, and so $\Omega$
is a complex algebraic variety.   We view $\Omega$ as
a component in the algebraic variety $\Omega(G)$.

We recall the {\it extended quotient}.   Let the finite group $\Gamma$
act on the space $X$.   Let $\widetilde{X} =
\{(x,\gamma) : \gamma x = x\}$, let $\Gamma$ act on $\widehat{X}$ by
$\gamma_1(x,\gamma) = (\gamma_1x, \gamma_1 \gamma \gamma_1^{-1})$. Then
$\widetilde{X}/\Gamma$ is the  extended quotient of $X$ by $\Gamma$.
There is a canonical projection $\widetilde{X}/\Gamma \rightarrow X/\Gamma$.

The Bernstein variety $\Omega(G)$ is the disjoint union
of ordinary quotients. We now
replace the ordinary quotient by the extended quotient to create
a new variety $\Omega^+(G)$. So we have
$$
\Omega(G) = \bigsqcup D/W(M, D)\; \; \mbox{\rm and}
\; \;   \Omega^+(G) = \bigsqcup \widetilde{D}/W(M, D)
$$

Let $\Omega$ be a component in the Bernstein variety
$\Omega(GL(n))$, and let $\mathcal{H}(G) = \bigoplus
\mathcal{H}(\Omega)$ be the Bernstein decomposition of the Hecke
algebra.

Let
\[
\Pi(\Omega) = (inf.ch.)^{-1}\Omega.
\]

Then $\Pi(\Omega)$ is a smooth complex algebraic variety with
finitely many irreducible components.   We have the following
Bernstein decomposition of $\Pi(G)$:
\[
\Pi(G) = \bigsqcup \Pi(\Omega).
\]

Let $M$ be a compact $C^{\infty}$ manifold.   Then $C^{\infty}(M)$
is a Fr\'{e}chet algebra, and we have Connes' fundamental theorem
\cite[Theorem 2, p. 208]{C}:
\[
\HP_*(C^{\infty}(M)) \cong \RH^*(M ; \mathbb{C}).
\]

Now the ideal $\mathcal{H}(\Omega)$ is a purely algebraic object, and,
in computing its periodic cyclic homology, we would hope to find an
algebraic variety to play the role of the manifold $M$.
This algebraic variety is $\Pi(\Omega)$.

\begin{theorem}\label{2.1}
   Let $\Omega$ be a component in the Bernstein variety
$\Omega(G)$.   Then the periodic cyclic homology of
$\mathcal{H}(G)$ is isomorphic to the periodised de Rham cohomology
of $\Pi(\Omega)$:

$$ \HP_* \,\mathcal{H}(\Omega) \cong \RH^*(\Pi(\Omega); \mathbb{C}). $$
\end{theorem}

\begin{proof}
We can think of $\Omega$ as a vector
$({\tau}_1, \ldots, {\tau}_r)$ of irreducible supercuspidal
representations of smaller general linear groups, the entries of this
vector being only determined up to tensoring with unramified
quasicharacters and permutation.   If the vector is equivalent to
$(\sigma_1, \ldots, \sigma_1, \ldots, \sigma_r, \ldots, \sigma_r)$
with $\sigma_j$ repeated $e_j$ times, $1\leq j\leq r$, and
$\sigma_1, \ldots, \sigma_r$ are pairwise distinct, then
we say that $\Omega$ has \emph{exponents} $e_1, \ldots, e_r$.

  Then there is a Morita equivalence
\[\mathcal{H}(\Omega) \sim \mathcal{H}(e_1,q_1) \otimes \ldots \otimes
\mathcal{H}(e_r,q_r)\]
where $q_1, \ldots, q_r$ are natural number invariants
attached to $\Omega$.

This result is due to Bushnell-Kutzko
\cite{BK1, BK2, BK3}.   We describe the steps
in the proof. Let $(\rho, W)$ be an
irreducible smooth representation of the compact open subgroup
$K$ of $G$.   As in \cite[4.2]{BK2}, the pair $(K,\rho)$ is an $\Omega$-type
in $G$ if and only if, for $(\pi, V) \in \Pi(G)$, we have
$\text{inf.ch.}(\pi) \in \Omega$ if and only if $\pi$ contains $\rho$.
The existence of an $\Omega$-type in $GL(n)$, for each component
$\Omega$ in $\Omega(GL(n))$, is established in
\cite[1.1]{BK3}.   So let
$(K, \rho)$ be an $\Omega$-type in $GL(n)$.   As in \cite[2.9]{BK2}, let

\[
e_{\rho}(x) = (\text{vol} K)^{-1}(\dim\,\rho)\,\tr_W(\rho(x^{-1}))
\]
for $x\in K$ and $0$ otherwise.

Then $e_{\rho}$ is an idempotent in the Hecke algebra $\mathcal{H}(G)$.
Then we have
\[\mathcal{H}(\Omega) \cong \mathcal{H}(G) * e_{\rho} * \mathcal{H}(G)\]
as in \cite[4.3]{BK2} and the two-sided ideal
$\mathcal{H}(G) * e_{\rho} * \mathcal{H}(G)$ is Morita
equivalent to $e_{\rho} * \mathcal{H}(G) * e_{\rho}$.
Now let $\mathcal{H}(K, \rho)$ be the endomorphism-valued Hecke algebra
attached to the semisimple type $(K, \rho)$.   By \cite[2.12]{BK2}
we have a canonical isomorphism of unital $\mathbb{C}$-algebras :

\[\mathcal{H}(G, \rho) \otimes_{\mathbb{C}} \End_{\mathbb{C}}W
\cong e_{\rho} * \mathcal{H}(G) * e_{\rho}\]

so that $e_{\rho} * \mathcal{H}(G) * e_{\rho}$ is Morita equivalent to
$\mathcal{H}(G, \rho)$.   Now we quote the main
theorem for semisimple types in $GL(n)$ \cite[1.5]{BK3}:  there is an
isomorphism of unital $\mathbb{C}$-algebras

\[\mathcal{H}(G, \rho) \cong \mathcal{H}(G_1, \rho_1) \otimes
\ldots \otimes \mathcal{H}(G_r, \rho_r)\]

The factors $\mathcal{H}(G_i, \rho_i)$ are (extended) affine Hecke algebras
whose structure is given explicitly in \cite[5.6.6]{BK1}.   This
structure is in terms of generators and relations \cite[5.4.6]{BK1}.
So let $\mathcal{H}(e,q)$ denote the affine Hecke algebra
associated to the affine Weyl group $\mathbb{Z}^e \rtimes S_e$.
Putting all this together we obtain a Morita equivalence

\[\mathcal{H}(\Omega) \sim \mathcal{H}(e_1,q_1) \otimes \ldots \otimes
\mathcal{H}(e_r,q_r)\]
The natural numbers $q_1, \ldots, q_r$ are specified in
\cite[5.6.6]{BK1}.   They are the cardinalities of the residue
fields of certain extension fields $E_1/F, \ldots, E_r/F$.

Using the K\"{u}nneth formula the calculation of $\HP_*(\mathcal{H}(\Omega))$
is reduced to that of the affine
Hecke algebra $\mathcal{H}(e,q)$.
When $q=1$ the Hecke algebra $\mathcal{H}(e,q)$ is isomorphic to the
complex group algebra of the affine Weyl group
$\mathbb{Z}^e \rtimes S_e$.
We now quote the main result in \cite{BN1,BN2}:
\[
\HP_* \mathcal{H}(e,q) \cong \HP_* \mathcal{H}(e,1)
\cong \HP_* \mathbb{C}[\mathbb{Z}^e \rtimes S_e].
\]
We have also \cite{BC}\cite{BN1}
$$
\HP_* \bbc [\bbz^e\rtimes S_e] \simeq \RH^*(\widetilde{\bbt^e}/S_e ; \bbc)
$$
where $\widetilde{\bbt^e}/S_e$ is the extended quotient of the torus
$\bbt^e$ by the symmetric group $S_e$.
Therefore we have
$$
\HP_*(\calh(e,q)) = \RH^*(\widetilde{\bbt^e}/S_e ; \bbc).
$$

If $\Omega$ has exponents $e_1, \ldots, e_r$ then
$e_1 + \ldots +e_r = d(\Omega) = \dim_\mathbb{C} \, \Omega$, and
$W(\Omega)$ is a product of symmetric groups:
\[W(\Omega) = S_{e_1} \times \ldots \times S_{e_r}\]

Form the semidirect product $\mathbb{Z}^{d(\Omega)} \rtimes W(\Omega)$.
Then we have
\[
\HP_*(\mathcal{H}(\Omega)) \cong \HP_*(\mathbb{C}[\mathbb{Z}^{d(\Omega)}
\rtimes W(\Omega)]).\]

We therefore have
\[ \HP_j(\mathcal{H}(\Omega))
\cong \oplus_l \RH^{j+2l}(\widetilde{\mathbb{T}^{d(\Omega)}}/W(\Omega); \mathbb{C})
\]
with $j = 0,1$.   By \cite[p. 217]{BP2} we have $\Pi(\Omega) \cong
\Omega^+$.   It now follows that
$$
\HP_* (\calh(\Omega)) = \RH^*(\Pi(\Omega); \mathbb{C}).
$$
\end{proof}

\begin{lemma}   Let $\Omega$ be a component in the variety
$\Omega(G)$.   Then we have the dimension formula

\[
\dim_\bbc \HP_* \calh(\Omega) = 2^{k(\phi_1)- 1} + \ldots +
2^{k(\phi_r) - 1}.
\]
\end{lemma}

\begin{proof}Let
\[
\Phi(\Omega) = (inf.ch. \circ \pi_F)^{-1}\Omega = (\pi_q
\circ \alpha)^{-1} \Omega.
\]

Then $\Phi(\Omega)$ is a disjoint union of orbits
\cite[p.217]{BP2} and we have
\[
\Phi(\Omega) = \mathcal{O}(\phi_1) \sqcup \ldots \sqcup
\mathcal{O}(\phi_r) \cong \Omega^+.
\]

By Theorem 2.1 we have
$$
\HP_*(\calh(\Omega)) \simeq \RH^*(\Omega^+; \bbc)
$$
and the lemma easily follows.

\end{proof}

Theorem 2.1, combined with the calculation in \cite{BP1}, now leads to the next result.

\begin{theorem}   The inclusion $\mathcal{H}(G) \longrightarrow
\mathcal{S}(G)$ induces an isomorphism at the level of periodic
cyclic homology:
$$
\HP_*(\calh(G)) \simeq \HP_*(\cals(G)).
$$
\end{theorem}

\section{The $q$-projection}

Let $\Omega$ be a component in the Bernstein variety. This component
is an ordinary quotient $D/\Gamma$. We now consider the extended quotient
$\widetilde{D}/\Gamma = \bigsqcup D^\gamma/Z_\gamma$, where $D$ is the complex
torus $\bbc^{\times m}$. Let $\gamma$ be a permutation of $n$ letters
with cycle type
$$
\gamma = (1\ldots  \alpha_1)\cdots (1\ldots \alpha_r)
$$
where $\alpha_1 + \dots +\alpha_r = m$. On the fixed set $D^\gamma$
the map $\pi_q$, by definition, sends the element
$
(z_1, \ldots , z_1, \ldots , z_r, \ldots , z_r)
$
where $z_j$ is repeated $\alpha_j$ times, $1\leq j\leq r$, to the element
$$
(q^{(\alpha_1 - 1)/2}z_1, \ldots , q^{(1 - \alpha_1)/2}z_1, \ldots ,
q^{(\alpha_r - 1)/2}z_r, \ldots ,
q^{(1 -\alpha_r)/2}z_r)
$$
 The map $\pi_q$ induces a map from $D^\gamma/Z_\gamma$
to $D/\Gamma$, and so a map, still denoted $\pi_q$,
from the extended quotient $\widetilde{D}/\Gamma$ to the
ordinary quotient $D/\Gamma$. This creates a map $\pi_q$
from the extended Bernstein variety to the Bernstein variety:
$$
\pi_q:
\Omega^+(G) \longrightarrow \Omega(G).
$$

\begin{definition}
The map $\pi_q$ is called the \emph{$q$-projection}.
\end{definition}

The $q$-projection $\pi_q$ occurs in the following commutative diagram
\cite{BP2}:
$$
\begin{CD}
\Phi(G) @>>> \Pi (G) \\
@V{\alpha}VV @VV{\text{inf. ch.}}V \\
\Omega^+(G) @>{\pi_q}>> \Omega(G)
\end{CD}
$$

Let $A, B$ be commutative rings with $A
\subset B, 1 \in A$.  Then the element $x \in B$ is {\it integral}
over $A$ if there exist $a_1, \ldots, a_n \in A$ such that

$$x^n + a_1x^{n-1} + \ldots + a_n = 0.$$

Then $B$ is {\it integral} over $A$ if each $x \in B$ is integral
over $A$.   Let $X, Y$ be affine varieties, $f : X \longrightarrow
Y$ a regular map such that $f(X)$ is dense in $Y$.   Then the
pull-back $f^{\#}$ defines an isomorphic inclusion $\bbc[Y]
\longrightarrow \bbc[X]$.   We view $\bbc[Y]$ as a subring of
$\bbc[X]$ by means of $f^{\#}$.   Then $f$ is a {\it finite} map
if $\bbc[X]$ is integral over $\bbc[Y]$, see \cite{Sh}.
This implies that the pre-image $F^{-1}(y)$ of each point $y \in Y$ is
a finite set, and that, as $y$ moves in $Y$, the points in
$F^{-1}(y)$ may merge together but not disappear.    The map
$\bba^1 - \{0\} \longrightarrow \bba^1$ is the classic example of
a map which is {\it not} finite.

\begin{lemma}
Let X be a component in the extended variety $\Omega^+(G)$.   Then the  $q$-projection
$\pi_q$ is a finite map from $X$ onto its image $\pi_q(X)$.
\end{lemma}
\begin{proof}

Note that the fixed-point set $D^{\gamma}$ is a complex torus of dimension $r$,
that $\pi_q(D^{\gamma})$ is a torus of dimension $r$  and
that we have an isomorphism of affine varieties
$D^{\gamma} \cong \pi_q(D^{\gamma})$.   Let $X = D^{\gamma}/Z_{\gamma}, Y =
\pi_q(D^{\gamma})/\Gamma$ where $Z_{\gamma}$ is the $\Gamma$-centralizer of $\gamma$.
Now
each of $X$ and $Y$ is a quotient of the variety $D^{\gamma}$ by a finite group,
hence $X, Y$ are affine varieties
\cite[p.31]{Sh}.   We have $D^{\gamma} \longrightarrow X \longrightarrow Y$
and $\bbc[Y] \longrightarrow \bbc[X] \longrightarrow \bbc[D^{\gamma}]$.
According to \cite[p.61]{Sh}, $\bbc[D^{\gamma}]$ is integral over $\bbc[Y]$ since $Y =
D^{\gamma}/\Gamma$.   Therefore the subring $\bbc[X]$ is integral over
$\bbc[Y]$.   So the map $\pi_q : X \longrightarrow Y$ is finite.
\end{proof}

{\sc Example}. $GL(2)$.   Let $T$ be the diagonal subgroup of $G=GL(2)$ and let
$\Omega$ be the component in $\Omega \, G$ containing
the cuspidal pair
$(T,1)$.     Then $\sigma \in \Pi (GL(2))$ is {\it arithmetically
unramified} if $inf.ch.\sigma \in \Omega$.   If $\pi_F(\phi) =
\sigma$ then $\phi$ is a $2$-dimensional representation of
$\mathcal{L}_F$ and there are two possibilities:

{\it $\phi$ is reducible}, $\phi = \psi_1 \oplus \psi_2$ with
$\psi_1, \psi_2$ unramified quasicharacters of $W_F$.   So
$\psi_j(w) = z_j^{d(w)}, z_j \in \mathbb{C}^{\times}, j =1,2$.
We have $\pi_F(\phi) = Q(\psi_1, \psi_2)$ where $\psi_1$ does not precede
$\psi_2$.
  In particular we
obtain the
$1$-dimensional representations of $G$ as follows:
\[
\pi_F(| \;|^{1/2}\psi \oplus | \;|^{-1/2}\psi) = Q(| \;|^{1/2}\psi,
| \;|^{-1/2}\psi) = \psi \circ det.
\]
{\it $\phi$ is irreducible}, $\phi = \psi \otimes \text{spin}(1/2)$.
Then $\pi_F(\phi) = Q(\Delta)$ with $\Delta = \{|\;|^{-1/2} \psi, | \;|^{1/2} \psi\}$
so $\pi_F(\phi) = \psi \otimes St(2)$
where $St(2)$ is the Steinberg representation of
$GL(2)$.

The orbit of $(T,1)$ is $D = (\mathbb{C}^{\times })^2$, and $W(T,D) =
\mathbb{Z}/2\mathbb{Z}$.   Then $\Omega \cong
(\mathbb{C}^{\times })^2/\mathbb{Z}/2\mathbb{Z}\cong Sym^2 \,
\mathbb{C}^\times$.   The extended quotient is $\Omega^{+} = Sym^2 \; \mathbb{C}
^{\times} \sqcup \mathbb{C}^{\times}$.
The $q-$projection works as follows:
\[
\pi_q : \{z_1, z_2\} \mapsto \{z_1, z_2\}
\]

\[
\pi_q : z \mapsto \{q^{1/2}z, q^{-1/2}z\}
\]
where $q$ is the cardinality of the residue field of $F$.

\mbox{}

{\sc Example}.  $GL(3)$.   In the above example, the $q$-projection is
stratified-injective, i.e. injective on each orbit type.   This is not so in general, as
shown by next example (due to J-F. Dat).   Let $T$ be the diagonal
subgroup of $GL(3)$ and let $\Omega$ be the component containing
the cuspidal pair $(T,1)$.   Then $\Omega = Sym^3 \, \mathbb{C}^{\times}$
and \[\Omega^+ = Sym^3 \, \mathbb{C}^{\times} \sqcup (\mathbb{C}^{\times})^2 \sqcup
\mathbb{C}^{\times}\].

The map $\pi_q$ works as follows:
\[
\{z_1, z_2, z_3\} \mapsto \{z_1, z_2, z_3\}
\]
\[(z, w, w) \mapsto \{z,
q^{1/2}w, q^{-1/2}w\}
\]
\[(z, z, z) \mapsto \{qz, z, q^{-1}z\}.
\]

Consider the $L$-parameter
\[
\phi = \psi_1 \otimes 1 \oplus  \psi_2 \otimes \text{spin}(1/2) \in \Phi(GL(3)).
\]

If $\psi(w) = z^{d(w)}$ then we will write $\psi = z$.   With this
understood, let
\[
\phi_1 = q \otimes 1 \oplus q^{-1/2} \otimes \text{spin}(1/2)
\]
\[
\phi_2 = q^{-1} \otimes 1 \oplus q^{1/2} \otimes \text{spin}(1/2).
\]

Then $\alpha(\phi_1), \alpha(\phi_2)$ are distinct points in the same stratum of
the extended quotient, but their image under the $q$-projection $\pi_q$ is the single
point $\{q^{-1}, 1, q \} \in Sym^3 \, \mathbb{C}^{\times}$.

Let
\[
\phi_3 = 1 \otimes \text{spin}(3/2)
\]
\[
\phi_4 = q^{-1} \otimes 1 \oplus 1 \otimes 1 \oplus q \otimes 1.
\]

Then the distinct $L$-parameters $\phi_1, \phi_2, \phi_3, \phi_4$ all have
the same image under the $q$-projection $\pi_q$.

\section{The Geometric Counterpart}

Given an $L$-parameter $\phi: \call_F \ra GL(n, \bbc)$ we have
$$
\phi = \phi_1 \oplus \ldots \oplus \phi_m
$$
with each $\phi_j$ an irreducible representation. We have
$\phi_j = \rho_j \otimes \spin (j)$ where each $\rho_j$ is an irreducible
representation of the Weil group $W_F$. We shall assume that $\det \rho_j$
is a unitary character.
Let $\calo(\phi)$ be the orbit of $\phi$ as in Definition \ref{orbit}.
The map $\calo(\phi) \ra \calo ^t(\phi)$ is now defined as follows
$$
\psi_1\phi_1 \oplus \ldots \oplus \psi_m\phi_m \mapsto
|\psi_1|^{-1}\cdot \psi_1\phi_1\oplus \ldots \oplus
|\psi_m|^{-1}\cdot \psi_m \phi_m.
$$
This map is a deformation retraction of the complex orbit $\calo(\phi)$
onto the compact orbit $\calo^t(\phi)$. Since $\Pi(G)$ is a disjoint
union of such complex orbits this formula determines, via the local
Langlands correspondence for $GL(n)$,  a deformation retraction
of $\Pi(G)$ onto the tempered dual $\Pi^t(GL(n))$.

The results of the paper may be summarized as follows.

\begin{theorem}
The smooth dual $\Pi (GL(n))$ has a natural complex structure: it is a complex
 manifold, with infinitely
many components.  In the context
 of the commutative diagram
$$
\begin{CD}
K_*^{\text{top}}(G)  @>{\mu}>>    K_*(C^*_r(G)) \\
@V{\text{ch}}VV                   @VV{\text{ch}}V  \\
\HP_*(\mathcal{H}(G)) @>{\imath_*}>> \HP_*(\mathcal{S}(G)) \\
@VVV     @VVV \\
\RH^*_c(\Pi (G); \mathbb{C}) @>>> \RH^*_c(\Pi^t(G);
\mathbb{C})
\end{CD}
$$
the Baum-Connes  map has a
geometric counterpart: it is induced by  the deformation retraction of $\Pi
(GL(n))$ onto the tempered dual $\Pi^t(GL(n))$.
\end{theorem}

\section{Supercuspidal representations of $GL(n)$}

In this section we track the fate of supercuspidal representations
of $GL(n)$ through the commutative diagram.
Let $\rho$ be an irreducible $n$-dimensional complex representation
of the Weil group $W_F$ such that $\det\rho $ is a unitary character and
let $\phi = \rho\otimes 1$.
Then $\phi$ is the $L$-parameter for
a pre-unitary supercuspidal representation $\omega$ of $GL(n)$. Let
$\calo(\phi)$
be the orbit of $\phi$ and $\calo^t (\phi)$ be the compact orbit
of $\phi$. Then $\calo(\phi)$ is a component in
the Bernstein variety isomorphic to $\bbc^\times$ and $\calo^t (\phi)$
is a component in the tempered dual, isomorphic to $\bbt$.
The $L$-parameter $\phi$ now determines the following data.

Let $(J, \lambda)$ be a maximal simple type for $\omega$ in the sense of
Bushnell and Kutzko \cite[chapter 6]{BK1}. Then $J$ is a compact open subgroup
of $G$ and $\lambda $ is a smooth irreducible complex representation
of $J$.

We will write
\[
\mathbb{T} = \{\psi \otimes \omega : \psi \in \Psi^t(G)\}
\]
where $\Psi^t(G)$ denotes the group of unramified unitary characters of $G$.

\begin{theorem}
Let $K$ be a maximal compact subgroup of  $G$ containing $J$
and form the induced representation $W = \text{Ind}^K_J (\lambda)$.
We then have
$$
\ell^2(G\times _K W) \simeq
\text{Ind}^G_K(W) \simeq \text{Ind}^G_J(\lambda) \simeq \int_{\bbt}\pi
d\pi.
$$
\end{theorem}

\begin{proof}
The supercuspidal representation $\omega$ contains $\lambda$
and, modulo unramified unitary twist, is the only irreducible unitary
representation
with this property \cite[6.2.3]{BK1}.
Now the Ahn reciprocity theorem
expresses $\text{Ind}_J^G$ as a direct integral \cite[p.58]{Lipsman}:
\[
\text{Ind}_J^G(\lambda) = \int n(\pi, \lambda) \pi d\pi
\]
where $d\pi$ is Plancherel measure and $n(\pi, \lambda)$ is the multiplicity of
$\lambda$ in $\pi|_J$.   But the Hecke algebra of a maximal simple type is
commutative (a Laurent polynomial ring).   Therefore $\omega|_J$
contains $\lambda$ with multiplicity $1$ (thanks to C. Bushnell for
this remark).   We then have
$n(\psi \otimes \omega, \lambda) = 1$ for all $\psi \in \Psi^t(G)$.   We note that
Plancherel measure induces Haar measure on $\mathbb{T}$, see \cite{P}.

The affine building of $G$ is defined
as follows \cite[p. 49]{Tits}:
\[
\beta G = \mathbb{R} \times \beta SL(n)
\]
where $g \in G$ acts on the affine line $\mathbb{R}$ via $t \mapsto
t + val(det(g))$.   Let $G^{\circ} = \{g \in G : val(det(g)) = 0 \}$.   We use the standard
model for $\beta SL(n)$ in terms of equivalence classes of
$\mathbf{o}_F$ - lattices in the $n$-dimensional $F$-vector space $V$.   Then the
vertices of $\beta SL(n)$ are in bijection with the maximal
compact subgroups of $G^{\circ}$, see \cite[9.3]{Ronan}.   Let $P \in \beta G$ be the vertex
for which the isotropy subgroup is $K = GL(n, \mathbf{o}_F)$.   Then the
$G$-orbit of $P$ is the set of all vertices in $\beta G$ and
the discrete space $G/K$ can be identified
with the set of vertices in the affine building $\beta G$.
Now the base space of the associated vector
bundle $G\times_KW$ is the discrete coset space $G/K$, and the
Hilbert space of $\ell^2$-sections of this homogeneous vector
bundle is a realization of the induced representation
$\text{Ind}^G_K(W)$.
\end{proof}

The $C_0(\beta G)$-module structure is defined as follows.
Let $f \in C_0(\beta G)$, $s \in \ell^2(G \times_KW)$ and define
\[
(fs)(v) = f(v)s(v)
\]
for each vertex $v \in \beta G$.  We proceed to construct a $K$-cycle in degree $0$.
This $K$-cycle is \[(C_0(\beta G), \ell ^2 (G\times_K W) \oplus 0, 0)\]
interpreted as a $\bbz/2\bbz$-graded module.
This triple satisfies
the properties of a (pre)-Fredholm module \cite[IV]{C}
 and so creates
an element in $K_0^{\text{top}}(G)$. By Theorem 5.1  this generator
creates a free
$C(\bbt)$-module of rank 1, and so provides a generator in $K_0(C^*_r(G))$.

{\bf 5.2}   The Hecke algebra of the maximal simple type $(J, \lambda)$ is commutative
(the Laurent polynomials in one complex variable).  The periodic cyclic homology of this algebra
is generated by $1$ in degree zero and $dz/z$ in degree $1$.

The corresponding summand of the Schwartz algebra $\cals(G)$ is Morita equivalent to
the Fr\'{e}chet algebra $C^\infty(\bbt)$.   By an elementary
application  of Connes' theorem \cite[ Theorem 2, p. 208]{C},
the periodic cyclic homology of this Fr\'{e}chet algebra is generated
by $1$ in degree $0$ and $d\theta$ in degree $1$.

{\bf 5.3}   The corresponding component in the Bernstein variety is a copy
of $\bbc^\times$. The cohomology of $\bbc^\times$ is generated
by $1$ in degree $0$ and $d\theta$ in degree $1$.

The corresponding component in the tempered dual is the
circle $\mathbb{T}$.  The cohomology of $\mathbb{T}$ is generated by $1$ in
degree $0$ and $d\theta $ in degree $1$.


\begin{thebibliography}{99}
\bibitem{BHP} P.~Baum, N.~Higson, R.~J.~Plymen, A proof of the Baum-Connes
conjecture
for $p$-adic $GL(n)$, C.~R. Acad. Sci. Paris 325 (1997) 171-176.
\bibitem{BHP2} P.~Baum, N.~Higson, R.~J.~Plymen, Representation
theory of $p$-adic groups: a view from operator algebras,
Proc. Symp. Pure Math. 68 (2000) 111 -- 149.
\bibitem{BC}  P. Baum, A. Connes, Chern character for discrete
groups, in: A F\^{e}te of Topology, Academic Press, New York, 1988, p.163 --232.
\bibitem{BN1}  P. Baum, V. Nistor, The periodic cyclic homology of
Iwahori-Hecke algebras, C. R. Acad. Sci. Paris 332 (2001) 1 --6.
\bibitem{BN2}  P. Baum, V. Nistor, The periodic cyclic homology of
Iwahori-Hecke algebras, preprint 2001.
\bibitem{Be} J.~Bernstein, Representations of $p$-adic groups, Notes by
K.E. Rumelhart, Harvard University, 1992.
\bibitem{Borel} A. Borel, Automorphic $L$-functions, Proc. Symp. Pure
Math. 33, part 2, (1979) 27 -- 61.
\bibitem{BP1}  J. Brodzki, R.J. Plymen,  Periodic cyclic
homology of certain nuclear algebras, C. R. Acad. Sci. Paris., 329 (1999),
671--676.
\bibitem{BP2} J.~Brodzki, R.~J.~Plymen, Geometry of the smooth dual of
$GL(n)$, C. R. Acad. Sci. Paris., 331 (2000), 213--218.
331(2000) 213 -- 218.
\bibitem{BP3}  J. Brodzki, R.J. Plymen, Chern character for the
$p$-adic group $GL(n)$, Bull. London. Math. Soc. 34 (2002), to
appear.
\bibitem{BK1}  C.J. Bushnell and P.C. Kutzko, The admissible dual
of $GL(N)$ via compact open subgroups, Annals of Math. Studies
129, Princeton University Press, Princeton, 1993.
\bibitem{BK2}  C.J. Bushnell and P.C. Kutzko, Smooth
representations of $p$-adic reductive groups: Structure theory via
types, Proc. London Math. Soc 77(1998) 582 -- 634.
\bibitem{BK3}  C.J. Bushnell and P.C. Kutzko, Semisimple types for
$GL(N)$, Compositio Math. 119 (1999) 53 --97.
\bibitem{C} A.~Connes, Noncommutative Geometry, Academic Press,
New York 1994.
\bibitem{HT} M.~Harris, R.~Taylor, On the geometry and cohomology
of some simple Shimura varieties, Ann. Math. Study 151, Princeton
University Press, to appear.
\bibitem{Hen} G.~Henniart, Une preuve simple des conjectures de Langlands
pour $GL(n)$ sur un corps $p$-adique,
Invent. Math. 139 (2000) 439 --455.
\bibitem{Knapp} A.~W.~Knapp, Introduction to the Langlands program,
Proc. Symp. Pure Math. 61 (1997) 245--302.
\bibitem{Ku}  S.~S.~Kudla, The local Langlands correspondence,
Proc. Symp. Pure Math. 55 (1994) 365 --391.
\bibitem{L} V.~Lafforgue, Une d\'{e}monstration de la conjecture
de Baum-Connes pour les groupes r\'{e}ductifs sur un corps $p$-adique
et pour certains groupes discrets poss\'{e}dant la prori\'{e}t\'{e}
(T), C. R. Acad. Sci. Paris 327 (1998), 439--444.
\bibitem{Stuhler} G. Laumon, M. Rapoport, U. Stuhler,
$\mathcal{D}$-elliptic sheaves and the Langlands correspondence,
Invent. Math. 113 (1993) 217 --338.
\bibitem{Lipsman} R.L. Lipsman, Group representations, Lecture Notes in Math. 388, Springer, 1974.
\bibitem{N} J.~Neukirch, Algebraic Number Theory, Springer, Berlin, 1999.
\bibitem{N2} J. Neukirch, Local class field theory, Springer, Berlin 1986.
\bibitem{P}  R.J. Plymen, Reduced $C^*$-algebra of the $p$-adic
group $GL(n) \, II$, arXiv:math.OA/0110018.
\bibitem{Ronan} M. Ronan, Lectures on buildings, Perspectives in
Math. 7, Academic Press, 1989.
\bibitem{Sh}  I.~R.~Shafarevich, Basic algebraic geometry 1, Springer,
Berlin, 1994.
\bibitem{Sm} K.~E.~Smith, L. Kahanp\"{a}\"{a}, P. Kek\"{a}l\"{a}inen, W. Traves, An invitation to algebraic geometry,
Universitext, Springer, Berlin, 2000.
\bibitem{T} J.~Tate, Number theoretic background, Proc. Symp. Pure
Math. 33 (1979) part 2, 3 - 26.
\bibitem{Tits} J. Tits, Reductive groups over local fields, Proc.
Symp. Pure Math. 33 (1979) part 1, 29 -- 69.
\end{thebibliography}
\end{document}